\documentclass{amsart}
\usepackage{amsmath, amssymb, latexsym, amscd}
\usepackage{graphicx, psfrag}

\newcommand{\bC}{\mathbb{C}}
\newcommand{\bE}{\mathbb{E}}
\newcommand{\bL}{\mathbb{L}}
\newcommand{\bP}{\mathbb{P}}
\newcommand{\bQ}{\mathbb{Q}}
\newcommand{\bR}{\mathbb{R}}
\newcommand{\bZ}{\mathbb{Z}}

\newcommand{\cM}{\mathcal{M}}
\newcommand{\cO}{\mathcal{O}}

\newcommand{\cW}{\mathcal{W}}

\newcommand{\Mbar}{\overline{\mathcal{M}} }

\newcommand{\Aut}{\mathrm{Aut}}

\newcommand{\Ker}{\mathrm{Ker}}
\newcommand{\vir}{\mathrm{vir}}
\newcommand{\rank}{\mathrm{rank}}

\newcommand{\lam}{\lambda}

\newcommand{\vd}{\vec{d}}
\newcommand{\vnu}{{\vec{\nu}}}
\newcommand{\vmu}{{\vec{\mu}}}
\newcommand{\vn}{\mathbf{n}}

\newcommand{\bu}{\bullet}
\newcommand{\tF}{\tilde{F}}
\newcommand{\tC}{\tilde{C}} 

\newcommand{\up}[1]{ {{#1}^1,{#1}^2,{#1}^3} }

\newtheorem{thm}{Theorem}

\newtheorem{lm}[thm]{Lemma}

\newtheorem{pro}[thm]{Proposition}
\newtheorem{ex}[thm]{Example}

\begin{document}

\title{Gromov-Witten Invariants  of Toric Calabi-Yau Threefolds}
\author{Chiu-Chu Melissa Liu}
\address{Department of Mathematics,
Columbia University, New York, NY 10027, USA}
\email{ccliu@math.columbia.edu}

\subjclass[2000]{14N35}
\keywords{topological vertex, Gromov-Witten invariants, Calabi-Yau threefolds}

\begin{abstract}
Based on the large $N$ duality relating topological string theory on Calabi-Yau 3-folds
and Chern-Simons theory on 3-manifolds, M. Aganagic, A. Klemm, M. Mari\~{n}o and C. Vafa
proposed the topological vertex, an algorithm on computing Gromov-Witten invariants
in all genera of any non-singular toric Calabi-Yau 3-fold.  
In this expository article, we describe the mathematical theory of
the topological vertex developed by J. Li, K. Liu, J. Zhou,
and the author. 
\end{abstract}

\maketitle

\tableofcontents

\section{Gromov-Witten invariants of Calabi-Yau 3-folds}

\subsection{Symplectic and algebraic Gromov-Witten invariants}
We start with a general setting. 
Let $(X,\omega, J)$ be a compact K\"{a}hler manifold, 
where $\omega$ is the K\"{a}hler form, and $J$ is 
the complex structure. 
For our purpose, we may assume that $X$ is a projective 
manifold, i.e., a compact complex submanifold of 
some complex projective space $\bP^m$, and $\omega$
is the restriction of the Fubini-Study K\"{a}hler
form on $\bP^m$. In this case, one may work in
complex algebraic geometry.

Intuitively, symplectic Gromov-Witten invariants count parametrized
holomorphic curves in $X$. When $X$ is projective, algebraic Gromov-Witten 
invariants count parametrized complex algebraic curves in $X$,
and should coincide with symplectic Gromov-Witten invariants.

\subsection{Moduli space of stable maps}
Gromov-Witten invariants can be viewed as intersection numbers
on moduli spaces of parametrized holomorphic (complex algebraic) curves in $X$.
Let
$\cM_{g,0}(X,\vd)$ be  the moduli space of holomorphic
maps (morphisms) $f:C\to X$, where $C$ is a compact Riemann surface 
(smooth complex algebraic curve) of
genus $g$, and $f_*[C]=\vd\in H_2(X;\bZ)$. 
We call $\vd$ the {\em degree} of the map. Two maps are 
equivalent if they differ by an automorphism of the domain $C$.

To do intersection theory, we should compactify
$\cM_{g,0}(X,\vd)$. The standard compactification in  Gromov-Witten theory is 
$\Mbar_{g,0}(X,\vd)$, 
the Kontsevich's moduli space of stable maps
$f:C\to X$ of genus $g$, degree $\vd$ \cite{Ko2}, where the domain curve
$C$ has at most nodal singularities, and 
the map $f$ is {\em stable} in the sense that the automorphism group
of $f$ is finite.  When $X$ is projective, the compactified
moduli space $\Mbar_{g,0}(X.\vd)$ is a {\em proper Deligne-Mumford
stack} (in algebraic geometry),
or a compact, Hausdorff, singular orbifold (in differential geometry). 
(See \cite{Be-Ma, Fu-Pa}.)
Roughly,  ``proper'' corresponds to ``compact and Hausdorff'',
and ``Deligne-Mumford stack'' corresponds to 
``singular orbifold''. In our context, a singular orbifold
is a space which is locally of the form
$V/\Gamma$, where $V$ is the zero locus of polynomials
defined on an open set in $\bC^N$, and $\Gamma$ is 
a finite group acting on $V$. 

The moduli space  $\Mbar_{g,0}(X,\vd)$ is not a smooth manifold, so
it does not have a tangent bundle. However, it
has a {\em virtual tangent bundle}  which 
is the difference $E_0-E_1$ of two complex vector bundles
$E_0$ and $E_1$ over $\Mbar_{g,0}(X,\vd)$.
The virtual (complex) dimension of
$\Mbar_{g,0}(X,\vd)$ is defined to be
the rank of the virtual tangent bundle:
\begin{equation}\label{vdim}
\textup{vir. dim. }=\rank_\bC E_0 -\rank_\bC E_1=
\int_{\vd} c_1(T_X) +(\dim_\bC X-3)(1-g).
\end{equation}
The structure of a virtual tangent bundle 
in symplectic Gromov-Witten theory
corresponds to the structure of a {\em perfect
obstruction theory} in algebraic Gromov-Witten theory.
The rank of the virtual tangent bundle
is the virtual dimension of the
perfect obstruction theory.

\subsection{Gromov-Witten invariants of compact Calabi-Yau 3-folds}
When $X$ is a Calabi-Yau $n$-fold, in the
sense that $K_X=\Lambda^n T^*X$ is 
a trivial holomorphic line bundle over $X$,
we have $c_1(T_X)=0$. By the formula
\eqref{vdim}, the virtual dimension 
of $\Mbar_{g,0}(X,\vd)$ is $(n-3)(1-g)$, which is independent
of the degree $\vd$. In particular, when 
$X$ is a Calabi-Yau 3-fold, the virtual
dimension  of $\Mbar_{g,0}(X,\vd)$
is zero for any genus $g$ and any degree $\vd$. 
In this case, there is a {\em virtual fundamental class}
$$
[\Mbar_{g,0}(X,\vd)]^\vir\in H_0(\Mbar_{g,0}(X,\vd);\bQ).
$$ 
The virtual fundamental class 
has been constructed in a much more general setting
by Li-Tian \cite{Li-Ti1},  Behrend-Fantechi \cite{Be-Fa} 
in algebraic Gromov-Witten theory, and by Li-Tian \cite{Li-Ti2},
Fukaya-Ono \cite{Fu-On}, Ruan \cite{Ru}, Siebert \cite{Si}
(more recently, Hofer-Wysocki-Zehnder \cite{HWZ1, HWZ2, HWZ3}) in 
symplectic Gromov-Witten theory.

The genus $g$, degree $\vd$ {\em Gromov-Witten invariant} of 
a Calabi-Yau 3-fold $X$ is defined by 
\begin{equation}\label{eqn:Ngd}
N^X_{g,\vd}=\int_{[\Mbar_{g,0}(X,\vd)]^\vir} 1.
\end{equation}
where $\int$ stands for the pairing between
$H_0(\Mbar_{g,0}(X,\vd);\bQ)$ and 
$H^0(\Mbar_{g,0}(X,\vd);\bQ)$. 
If $\Mbar_{g,0}(X,\vd)$ were a compact complex
manifold of dimension zero, it would consist of
finitely many points, and the right hand side
of \eqref{eqn:Ngd} would be the number of points
in $\Mbar_{g,0}(X,\vd)$. In general,
$\Mbar_{g,0}(X,\vd)$ can be singular and
can have positive actual dimension. Then 
the right hand side of \eqref{eqn:Ngd} 
defines the ``virtual number'' of points
in $\Mbar_{g,0}(X,\vd)$. It is usually a rational
number instead of an integer because $\Mbar_{g,0}(X,\vd)$ is
an orbifold. For example, if a map has an automorphism 
group of order $2$, we count it as one half of a map
instead of one map. This is how fractional numbers arise.

To summarize, Gromov-Witten invariants are defined
for any smooth projective  Calabi-Yau 3-folds,
or more generally, any compact K\"{a}hler Calabi-Yau 3-folds,
for any genus and any degree.

\subsection{Gromov-Witten invariants of noncompact Calabi-Yau 3-folds}
The construction of the virtual fundamental
class requires two properties of the moduli space
$\Mbar_{g,0}(X,\vd)$: compactness (properness)
and the structure of a virtual tangent bundle
(perfect obstruction theory). When $X$ is not compact,
the moduli space $\Mbar_{g,0}(X,\vd)$ is usually
noncompact, but still equipped with a virtual tangent bundle
(perfect obstruction theory). Therefore, if  $X$ is noncompact
but $\Mbar_{g,0}(X,\vd)$ is compact for
a particular genus $g$ and degree $\vd$,  then the
Gromov-Witten invariant $N^X_{g,\vd}$ is defined
for the particular genus $g$ and degree $\vd$.

\begin{ex} \label{Pone}
\textup{Let  $X$ be the total space of $\cO_{\bP^1}(-1)\oplus\cO_{\bP^1}(-1)$. 
Then $X$ is a non-compact Calabi-Yau 3-fold. We have
$H_2(X;\bZ)\cong H_2(\bP^1;\bZ) =\bZ[\bP^1]$, 
where $[\bP^1]$ is the class of the zero section.
Any nonconstant holomorphic map from a compact Riemann
surface to $X$ factors through the embedding
$i_0:\bP^1\to X$ by the zero section. Therefore,
when $d\neq 0$,  the following two moduli spaces
are identical as topological spaces (or Deligne-Mumford
stacks):
\begin{equation}\label{eqn:Pone}
\Mbar_{g,0}(X,d[\bP^1])=\Mbar_{g,0}(\bP^1,d[\bP^1]).
\end{equation}
(They both are empty when $d<0$.)
However, they have different virtual tangent bundles
(perfect obstruction theories). The virtual dimension
of $\Mbar_{g,0}(X,d[\bP^1])$ is $0$
while that of $\Mbar_{g,0}(\bP^1, d[\bP^1])$ is $2(d+g-1)$.
For any $d>0$, the right hand side of \eqref{eqn:Pone} is 
compact because $\bP^1$ is compact. So  
$N_{g,d[\bP^1]}^X$ is defined for any genus $g$ 
and any $d\neq 0$, and is $0$ when $d<0$.}
\end{ex}

\begin{ex} \label{Ptwo}
\textup{Let $X$ be the total space of $\cO_{\bP^2}(-3)$.
Then $X$ is a noncompact Calabi-Yau 3-fold.
We have $H_2(X;\bZ)\cong H_2(\bP^2;\bZ)=\bZ\ell$, 
where $\ell$ is the class of a projective line $\bP^1$ in the zero section
$\bP^2$. Any nonconstant holomorphic map from a compact
Riemann surface to $X$ must factor through the embedding 
$i_0:\bP^2\to X$ by the zero section. Therefore, when $d\neq 0$,
the following two moduli spaces are identical 
as topological spaces (or Deligne-Mumford stacks):
\begin{equation}\label{eqn:Ptwo}
\Mbar_{g,0}(X,d\ell)=\Mbar_{g,0}(\bP^2,d\ell).
\end{equation}
(They both are empty when $d<0$.)
However, they have different virtual tangent bundles
(perfect obstruction theories).  The virtual dimension of
$\Mbar_{g,0}(X,d\ell)$ is $0$ while
that of $\Mbar_{g,0}(\bP^2,d\ell)$ is $3d+g-1$.
For any $d>0$,
the right hand side of \eqref{eqn:Ptwo} is compact
because $\bP^2$ is compact. So $N_{g,d\ell}^X$ is defined
for any genus $g$ and any $d\neq 0$, and is $0$ when
$d<0$.}
\end{ex}

Example \ref{Pone} and Example \ref{Ptwo} are examples of 
{\em toric} Calabi-Yau 3-folds. A Calabi-Yau 3-fold $X$ is 
{\em toric} if it contains
$(\bC^*)^3$ as an open dense subset, and the action
of $(\bC^*)^3$ on itself extends to $X$. In Example \ref{Pone},
removing $\infty$ from $\bP^1$, we obtain
a rank 2 holomorphic vector bundle over $\bP^1-\{\infty\} =\bC$, which must
be the trivial bundle, so the total space is $\bC\times \bC^2=\bC^3$. 
In Example \ref{Ptwo},
removing the line at infinity, we obtain  a holomorphic line bundle
over $\bP^2-\bP^1\cong \bC^2$, which must be
the trivial line bundle, so the total space is $\bC^2\times \bC=\bC^3$. 
In both examples, the inclusions $(\bC^*)^3\subset \bC^3\subset X$
are open and dense, and the $(\bC^*)^3$-action on itself extends to $X$.

\section{Traditional Algorithm in the Toric Case} \label{sec:tradition}

The traditional algorithm (the algorithm before
the ``topological vertex'') of computing 
Gromov-Witten invariants of toric Calabi-Yau
3-folds consists of two steps:
\begin{enumerate}
\item[T1.] {\em Localization} reduces Gromov-Witten invariants
of toric Calabi-Yau 3-folds to {\em Hodge integrals}, which 
are intersection numbers on moduli spaces of curves.
\item[T2.] {\em Hodge integrals} can be computed recursively.
\end{enumerate}
We will explain these two steps in this section.

\subsection{Localization}
The $(\bC^*)^3$-action on $X$ induces a $(\bC^*)^3$-action on
the moduli space $\Mbar_{g,0}(X,\vd)$ by moving the image
of a stable map.  Let $T$ be a subtorus of $(\bC^*)^3$
such that 
$$
\Mbar_{g,0}(X,\vd)^T  =\Mbar_{g,0}(X,\vd)^{(\bC^*)^3}
$$
where the left (resp. right) hand side is the set of $T$ 
(resp. $(\bC^*)^3$) fixed points in $\Mbar_{g,0}(X,\vd)$.
We have
\begin{equation}\label{eqn:localize}
N^X_{g,\vd}=\int_{[\Mbar_{g,0}(X,\vd)]^\vir} 1 
=\sum_F \int_{[F]^\vir}\frac{1}{e_T(N^\vir_F)}
\end{equation}
where the first equality is the definition
\eqref{eqn:Ngd}, and the second equality 
follows from the virtual localization formula
proved by Graber-Pandharipande \cite{Gr-Pa}. The sum is
over connected components $F$ of the fixed points
set $\Mbar_{g,0}(X,\vd)^T = \Mbar_{g,0}(X,\vd)^{(\bC^*)^3}$.
$[F]^\vir$ is the virtual fundamental class of
$F$ and $e_T(N^\vir_F)$ is the $T$-equivariant
Euler class of the virtual normal bundle $N^\vir_F$
of $F$ in $\Mbar_{g,0}(X,\vd)$. 
More explicitly, the virtual tangent bundle
of $\Mbar_{g,0}(X,\vd)$ is $T$-equivariant:
it  is of the form $E_0-E_1$,
where $E_0$ and $E_1$ are $T$-equivariant
complex vector bundles on $\Mbar_{g,0}(X,\vd)$.
For each fixed locus $F$, $T$ acts on $F$ trivially,
so $T$-action preserves the fibers of $E_i|_F$, the 
restriction of $E_i$ to $F$.
Let  $E_i^f$ and $E_i^m$ be the fixed and moving
parts of $E_i|_F$, respectively, so 
that $E_i|_F = E_i^f\oplus E_i^m$.
Then  $T_F^\vir = E_0^f - E_1^f$ is the virtual
tangent bundle of $F$, which defines 
the virtual fundamental class $[F]^\vir$,
and $N_F^\vir = E_0^m -E_1^m$
is the virtual normal bundle of $F$ in $\Mbar_{g,0}(X,\vd)$.

If $\Mbar_{g,0}(X,\vd)$
were a compact complex manifold, and each $F$ were
a compact complex submanifold, then $[F]^\vir$ 
would be  the usual fundamental class $[F]$ of $F$, $N^\vir_F$
would be  the usual normal bundle $N_F$ of $F$ in $\Mbar_{g,0}(X,\vd)$,
and the second equality in \eqref{eqn:localize} would be
the classical Atiyah-Bott localization formula \cite{At-Bo}.

In our case, $F$ is (up to finite cover) a product of 
moduli spaces of stable curves (see Section \ref{sec:hodge} below), 
which are smooth orbifolds (smooth Deligne-Mumford stacks), 
so it has a fundamental class $[F]\in H_*(F;\bQ)$. We have
$[F]=[F]^\vir$, and
\begin{equation}\label{eqn:E-E}
\int_{[F]^\vir}\frac{1}{e_T(N_F^\vir)}
=\int_{[F]}\frac{e_T(E_1^m)}{e_T(E_0^m)}.
\end{equation}
The integral on the right hand side of \eqref{eqn:E-E}
can be expressed in terms of Hodge integrals. The definition of Hodge integrals 
will be reviewed in the next subsection.

\subsection{Hodge integrals} \label{sec:hodge}

Let $\Mbar_{g,n}$ be the Deligne-Mumford compactification
of the moduli space of complex algebraic curves of genus $g$ with
$n$ marked points. 
A point in $\Mbar_{g,n}$ is represented by 
$(C,x_1,\ldots,x_n)$, where $C$ is a connected complex algebraic 
curve of arithmetic genus $g$ with at most nodal singularities,
$x_1,\ldots,x_n$ are distinct smooth points on $C$, 
and $(C,x_1,\ldots,x_n)$ is {\em stable} in the sense
that its automorphism group is finite. When $C$ is smooth,
it can be viewed as a connected compact Riemann surface of genus $g$. 
The moduli space $\Mbar_{g,n}$ is a proper, smooth Deligne-Mumford
stack (in algebraic geometry), or a compact, complex, smooth orbifold
(in differential geometry), of complex dimension $3g-3+n$. 
It is empty when $3g-3+n<0$.

The {\em Hodge bundle} $\bE$ is a rank $g$ vector bundle over
$\Mbar_{g,n}$ whose fiber over the moduli point $[(C,x_1,\ldots,x_n)]$
is $H^0(C,\omega_C)$, where $\omega_C$ is the dualizing
sheaf over $C$. When $C$ is smooth, $C$ can be viewed
as a compact Riemann surface
and $H^0(C,\omega_C)$ is the space of holomorphic
1-forms on $C$. The $\lambda$ classes are the Chern classes
of the Hodge bundle:
$$
\lambda_j=c_j(\bE)\in H^{2i}(\Mbar_{g,n};\bQ),\quad
j=1,\ldots,g.
$$
The cotangent line $T^*_{x_i}C$ of $C$ at the
$i$-th marked point $x_i$ gives rise to a line bundle
$\bL_i$ over $\Mbar_{g,n}$. 
The $\psi$ classes are the first Chern classes
of these line bundles:
$$
\psi_i=c_1(\bL_i)\in H^2(\Mbar_{g,n};\bQ),\quad i=1,\ldots,n.
$$
The $\lambda$ classes and $\psi$ classes lie in $H^*(\Mbar_{g,n};\bQ)$
instead of $H^*(\Mbar_{g,n};\bZ)$ because $\bE$ and $\bL_i$
are orbibundles on the orbifold $\Mbar_{g,n}$.

{\em Hodge integrals} are intersection numbers of $\lambda$ classes
and $\psi$ classes on $\Mbar_{g,n}$:
\begin{equation}\label{eqn:hodge}
\int_{\Mbar_{g,n}}\psi_1^{j_1}\cdots \psi_n^{j_n}
\lambda_1^{k_1}\cdots \lambda_g^{k_g} \in \bQ.
\end{equation}

The $\psi$ integrals (also known as {\em descendant integrals})
\begin{equation}\label{eqn:psi}
\int_{\Mbar_{g,n}}\psi_1^{j_1} \ldots \psi_n^{j_n}
\end{equation}
can be computed recursively by the Witten's conjecture \cite{Wi1}.
Witten's conjecture was  first proved by Kontsevich \cite{Ko1},
and there are now several alternative proofs 
\cite{Ok-Pa1, Mi, Ki-Liu, Ka-La, CLL}. 
Using Mumford's Grothendieck-Riemann-Roch calculations in \cite{Mu},
Faber showed in \cite{Fa} that general Hodge integrals \eqref{eqn:hodge}
can be uniquely reconstructed from descendant integrals \eqref{eqn:psi}.

\section{Physical Theory of the Topological Vertex}
\label{sec:AKMV}

Based on the {\em large $N$ duality} \cite{Go-Va1} between
the topological string theory on Calabi-Yau 3-folds
and the Chern-Simons theory on 3-manifolds,
Aganagic, Klemm, Mari\~{n}o, and Vafa proposed
the {\em topological vertex} \cite{AKMV}, an algorithm of
computing Gromov-Witten invariants in all genera of any 
smooth toric Calabi-Yau 3-folds.
Their algorithm can be summarized in the following
three steps.

\begin{enumerate}
\item[O1. ]{\em Topological vertex.}
There exist certain open Gromov-Witten invariants 
that count holomorphic maps from bordered Riemann surfaces
to $\bC^3$ with boundaries mapped to three Lagrangian
submanifolds $L_1, L_2, L_3$ (see Figure 1). 
Such invariants depend on the following discrete data:
\begin{enumerate}
\item[(i)] the topological
type of the domain, classified by the genus and  the number of
boundary circles;
\item[(ii)] the topological type of 
the map, described by a triple
of partitions $\vmu=(\up{\mu})$ where
$\mu^i=(\mu^i_1,\mu^i_2,\ldots)$ are degrees ("winding numbers") of 
boundary circles in $L_i\cong S^1\times \bC$;
\item[(iii)] the ``framing"
$n_i\in\bZ$ of the Lagrangian submanifolds $L_i$ ($i=1,2,3$).
\end{enumerate}
The {\em  topological vertex}
       $$C_{\vmu}(\lam;\vn)$$
       is a generating function of such invariants where
       one fixes the winding numbers $\vmu=(\up{\mu})$ and the framings
       $\vn=(n_1,n_2,n_3)$ and sums over the genus of the domain.
It can be viewed as local open Gromov-Witten invariants
of $D_1\cup D_2\cup D_3$ embedded in $(\bC^3, L_1\cup L_2\cup L_3)$
as in Figure 1.

\begin{figure}[h] \label{Figure1}
\begin{center}
\psfrag{L1}{\small $L_1\cong S^1\times \bC$}
\psfrag{L2}{\small $L_2$}
\psfrag{L3}{\small $L_3$}
\psfrag{D1}{\small $D_1$}
\psfrag{D2}{\small $D_2$}
\psfrag{D3}{\small $D_3$}
\psfrag{Z1}{\small $z_1$-axis}
\psfrag{Z2}{\small $z_2$-axis}
\psfrag{Z3}{\small $z_3$-axis}
\includegraphics[scale=0.8]{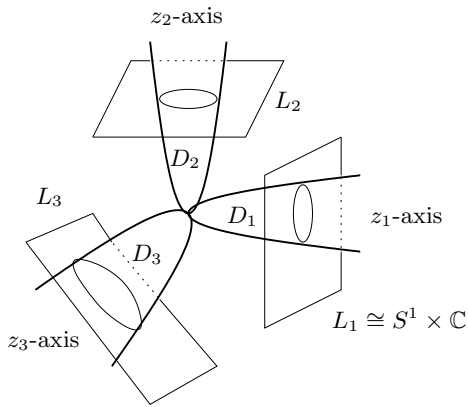}
\end{center}
\caption{$D_i$ is a holomorphic disk in $\bC^3$ with
boundary in the Lagrangian submanifold $L_i$.
$D_1=\{ (z_1,0,0)\mid z_1\in \bC, |z_1|\leq 1\}$, and
$L_1= \{ (\sqrt{1+|u|^2} e^{i\theta}, u, e^{-i\theta}\bar{u})\mid
e^{i\theta}\in S^1, u\in \bC\}$.
$D_2, D_3$ and $L_2, L_3$ can be obtained from
$D_1$ and $L_1$ by cyclic permutation of the three coordinates 
$z_1, z_2, z_3$.}
\end{figure}

\item[O2.] {\em Gluing algorithm.}
Any toric Calabi-Yau 3-fold $X$ can be constructed 
by gluing $\bC^3$ charts. 
The Gromov-Witten invariants of $X$
can be expressed in terms of local
open Gromov-Witten invariants
$C_{\vmu}(\lam;\vn)$ of $\bC^3$ by explicit gluing algorithm.

\item[O3.]{\em Closed formula.} 
By the large $N$ duality,
the topological vertex is given by
\begin{equation}\label{eqn:O3}
C_{\vmu}(\lam;\vn)=q^{\frac{1}{2}(\sum_{i=1}^3\kappa_{\mu^i}n_i)}
\cW_\vmu(q),\quad q=e^{\sqrt{-1}\lam},
\end{equation}
where $\kappa_\mu=\sum \mu_i(\mu_i-2i+1)$ 
for a partition $\mu=(\mu_1\geq \mu_2 \geq\cdots)$, 
and $\cW_\vmu(q)$ is a combinatorial expression related to 
the Chern-Simons link invariants. (Precise definition
of $\cW_\vmu(q)$ will be given in Section \ref{sec:W}.)
The left hand side of \eqref{eqn:O3} is an infinite series while
the right hand side of \eqref{eqn:O3} is a finite sum.
\end{enumerate}

The above algorithm is significantly more efficient than the
traditional algorithm described in Section \ref{sec:tradition}.
For example, to compute the degree 2 invariants of the
total space of $\cO_{\bP^2}(-3)$ using the traditional
algorithm, one computes one genus at a time; as
the genus increases, the computations soon become
too complicated to do by hand, so one needs
the aid of a computer. On the other hand,
using the algorithm of the topological vertex,
one can compute the generating function
of degree 2 invariants in all genera by hand.

Assuming O1 and the validity of open string virtual localization,
Diaconescu and Florea related $C_{\vmu}(\lambda;\vn)$ (at certain
fractional $n_i$) to Hodge integrals, and derived the gluing algorithms
in O2 by localization \cite{DF}.

When the toric Calabi-Yau threefold is the total space
of the canonical line bundle $K_S$ of a toric surfaces $S$ (e.g.
$K_{\bP^2}=\cO_{\bP^2}(-3)$), only
$\tC_{\mu^1,\mu^2,\emptyset}$ are required to evaluate their Gromov-Witten invariants. 
The algorithm in this case was described by Aganagic-Mari\~no-Vafa 
\cite{AMV}; an explicit formula was given by Iqbal \cite{Iq}
and derived by Zhou by localization, assuming a formula of 
Hodge integrals \cite{Zh1, Zh2}.

\section{Mathematical Theory of the Topological Vertex}
\label{sec:LLLZ}
J. Li, K. Liu, J. Zhou and the author developed
a mathematical theory of the topological vertex \cite{LLLZ} based
on relative Gromov-Witten theory. The relative Gromov-Witten
theory has been developed in symplectic geometry
by Li-Ruan \cite{Li-Ru} and Ionel-Parker \cite{Io-Pa1, Io-Pa2}. In our context, we need
to use the algebraic version developed by J. Li \cite{Li1, Li2}.
Our algorithm can be summarized as follows.

\begin{enumerate}
\item[R1.] We defined {\em formal relative Gromov-Witten invariants}
 for {\em relative formal toric Calabi-Yau (FTCY) 3-folds}.
 These invariants are refinements and generalizations of 
Gromov-Witten invariants of smooth toric Calabi-Yau 3-folds.

\item[R2.] Formal relative Gromov-Witten invariants satisfy the degeneration
formula. In particular, they can be expressed in terms of 
$\tilde{C}_{\vmu}(\lambda;\vn)$, formal relative Gromov-Witten
invariants of an indecomposable relative FTCY 3-fold. The degeneration
formula agrees with the gluing formula in O2, with $C_{\vmu}(\lambda;\vn)$
replaced by $\tilde{C}_\vmu(\lambda;\vn)$.

\item[R3.] $\tilde{C}_\vmu(\lambda;\vn)
=q^{(\sum_{i=1}^3 \kappa_{\mu^i} n_i)/2}\tilde{\cW}_{\vmu}(q)$,
where $\tilde{\cW}_{\vmu}(q)$ is a combinatorial expression
in terms of representations of symmetric groups. (The precise
definition of $\tilde{\cW}_\vmu(q)$ will be given in
Section \ref{sec:tW}.)
\end{enumerate}

We will describe this algorithm in detail in the remainder of this section.

\subsection{Locally planar trivalent graph} \label{sec:trivalent}
Let $X$ be a smooth toric Calabi-Yau 3-fold $X$.
For $k=0,1,2,3$, define $X^k$, the $k$-skeleton of $X$, 
to be the union of $k$-dimensional $(\bC^*)^3$ orbit closures.
Then
$$
X^0 \subset X^1 \subset X^2 \subset X^3=X
$$
where $X^0=X^{(\bC^*)^3}$, the set of $(\bC^*)^3$ fixed points
in $X$. 

When $X^0$ is nonempty, 
we may choose a distinguished
rank 2 subtorus $T$ of $(\bC^*)^3$ as follows. 
Pick any $(\bC^*)^3$ fixed point $p\in X^0$. Then
$(\bC^*)^3$ acts on $T_p X$, and its representation
on $\Lambda^3 T_p X\cong \bC$ gives a nontrivial 
irreducible character $\alpha:(\bC^*)^3 \to \bC^*$.
Define $T=\Ker(\alpha)\cong (\bC^*)^2$. Note that
the definition is independent of the choice of the
fixed point $p$ because $\Lambda^3 TX = K_X^{-1}$
is a trivial line bundle over $X$. 

$X^1$ is a configuration of rational curves. For example, 
let $X$ be  the total space of $\cO_{\bP^1}(-1)\oplus \cO_{\bP^1}(-1)$.
Then $X^1$ is the union of a projective line
$C_0\cong\bP^1$ and four complex lines $E_i \cong \bC$, 
$i=1, 2, 3, 4$, as shown on the left 
hand side of Figure 2. $C_0\cong \bP^1$ is the zero section, 
$p_0, p_1\in C_0$ are the two $T$-fixed points
on $X$, and $E_1, E_2$ (resp. $E_3, E_4$) are the
two $T$-invariant lines in the fiber
of $\cO_{\bP^1}(-1)\oplus \cO_{\bP^1}(-1)$ over
$p_0$ (resp. $p_1$).

We now associate to $X$ a locally planar trivalent
graph $\Gamma_X$. As an abstract graph, $\Gamma_X$ is 
determined by the configuration $X^1$:
each $T$-fixed point  corresponds to a vertex in $\Gamma_X$;
each $T$-invariant $\bP^1$ connecting two fixed points
corresponds to a compact edge (line segment) connecting
two vertices; each $T$-invariant $\bC$ containing
a fixed point corresponds to a noncompact edge (ray)
emanating from a vertex. Therefore all the vertices
in $\Gamma_X$ are trivalent. For example, the
graph for $\cO_{\bP^1}(-1)\oplus \cO_{\bP^1}(-1)$ is shown 
on the right hand side of Figure 2. The (local) embedding
of $\Gamma_X$ into $\bZ^2$ is determined by
the $T$-action on $X^1$: the slope of an edge
is determined by the $T$-action on the corresponding
irreducible component of $X^1$.
On the right hand side of Figure 2, 
the vectors $(1,0),(0,1),(-1,-1)$
(resp. $(-1,0), (0,-1), (1,1)$)
at the vertex $v_0$ (resp. $v_1$) correspond
to weights $t_1, t_2, t_1^{-1}t_2^{-1}$
(resp. $t_1^{-1}, t_2^{-1}, t_1 t_2$)
of the $T$-actions on $T_{p_0}C_0$, $T_{p_0}E_1$, 
$T_{p_0}E_2$ (resp. $T_{p_1}C_0$, $T_{p_1}E_3$, 
$T_{p_1}E_4$). The sum of  the three
vectors at $v_0$ (resp. $v_1$) is zero because
 $T$ acts trivially
on $\Lambda^3 T_{p_0}X$ (resp. $\Lambda^3 T_{p_1}X$). 

\begin{figure}[h] \label{Figure2}
\begin{center}
\psfrag{C0}{\small $C_0$}
\psfrag{E1}{\small $E_1$}
\psfrag{E2}{\small $E_2$}
\psfrag{E3}{\small $E_3$}
\psfrag{E4}{\small $E_4$}
\psfrag{p0}{\small $p_0$}
\psfrag{p1}{\small $p_1$}
\psfrag{v0}{\small $v_0$}
\psfrag{v1}{\small $v_1$}
\includegraphics[scale=0.6]{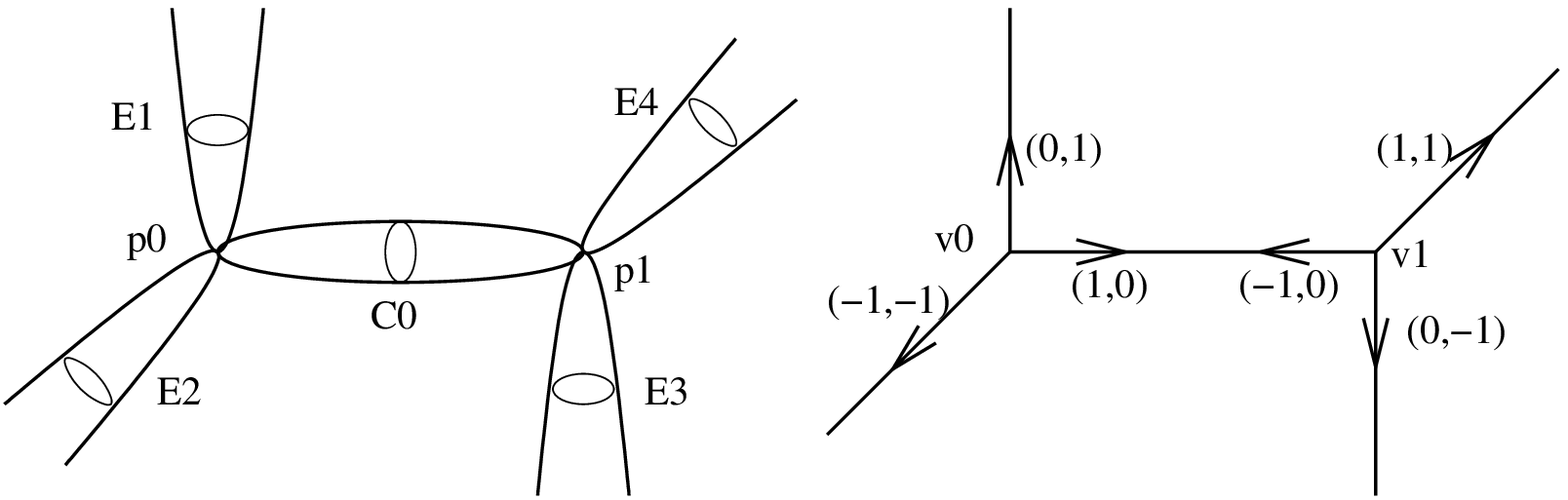}
\end{center}
\caption{$\cO_{\bP^1}(-1)\oplus \cO_{\bP^1}(-1)$.}
\end{figure}

More generally, given any toric Calabi-Yau 3-fold
$X$ and any $C_0\cong \bP^1$ which is an irreducible
component of $X^1$, the degree of the normal
bundle $N_{C_0/X}$ must be $-2$,
So $N_{C_0/X}\cong\cO_{\bP^1}(n)\oplus \cO_{\bP^1}(-n-2)$
for some $n\in \bZ$. We have a $T$-equivariant
open embedding $N_{C_0/X}\hookrightarrow X$. The
weights at one fixed point are $w_1, w_2, w_3 \in \bZ^2$,
where $w_1,w_2$ form a $\bZ$-basis of $\bZ^2$ and
$w_3=-w_1-w_2$. The weights at the other fixed point
is determined by $w_1,w_2, w_3$ and the degree $n$
(see Figure 3). Putting together the graphs of all $T$-invariant
$\bP^1$ in $X$, we obtain a locally planar trivalent
graph $\Gamma_X$. Some examples are shown in Figure 4.

\begin{figure}[h]  \label{Figure3}
\begin{center}
\psfrag{w1}{\small$w_1$}
\psfrag{w2}{\small$w_2$}
\psfrag{w3}{\small$w_3=-w_1-w_2$}
\psfrag{-w1}{\small$-w_1$}
\psfrag{w2-nw1}{\small$w_2-nw_1$}
\psfrag{w3+(n+2)w1}{\small$w_3+(n+2)w_1$} 
\includegraphics[scale=0.4]{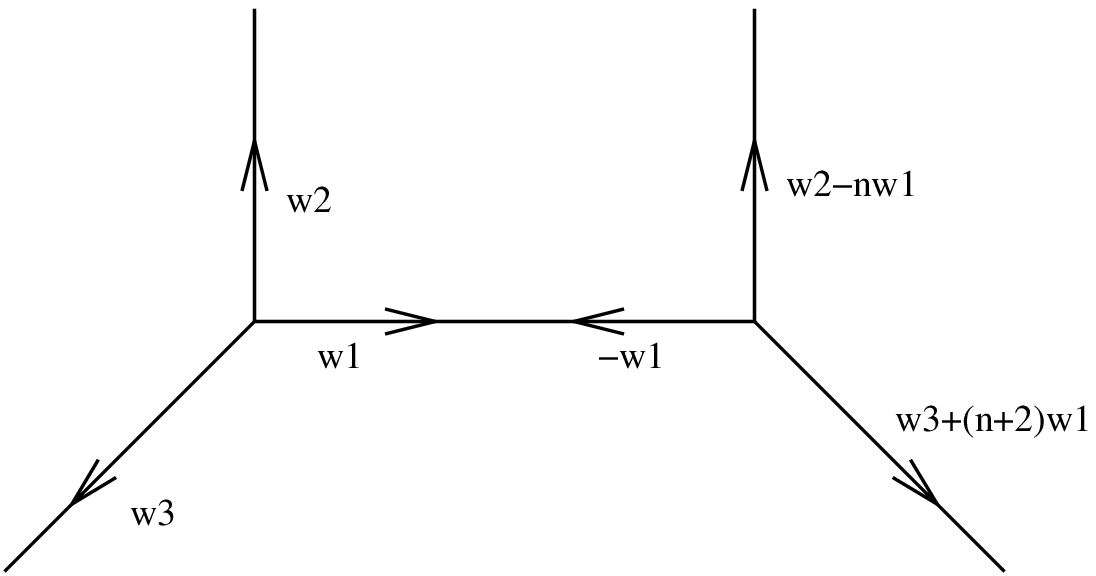}
\caption{$N_{C_0/X}=\cO_{\bP^1}(n)\oplus \cO_{\bP^1}(-n-2)$}
\end{center}
\end{figure}

\begin{figure}[h] \label{Figure4}
\begin{center}
\psfrag{P2}{$K_{\bP^2}=\cO_{\bP^2}(-3)$} 
\psfrag{P1P1}{$K_{\bP^1\times\bP^1}$}
\includegraphics[scale=0.5]{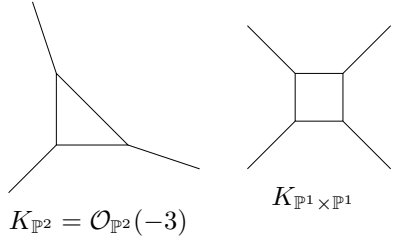}
\caption{Locally planar trivalent graphs}
\end{center}
\end{figure}

Conversely, from the graph in Figure 3 we can read off
the degrees of the two summands of 
$N_{C_0/X}$ together with the $T$-action
on $N_{C_0/X}$. Therefore, from $\Gamma_X$ we may recover
a $T$-equivariant formal neighborhood $\hat{X}^1$ of $X^1$ in $X$ 
(in algebraic geometry), or a $T_\bR$-equivariant
tubular neighborhood of $X^1$ in $X$ (in differential
geometry), where $T_\bR\cong U(1)^2$ is the maximal
compact subgroup of $T\cong (\bC^*)^2$.
The $T$-equivariant formal neighborhood $\hat{X}^1$
contains all the information needed to compute
all Gromov-Witten invariants $N_{g,\vd}^X$ of $X$
(using the traditional algorithm described in
Section \ref{sec:tradition}), 
because a point in $\Mbar_{g,0}(X,\vd)^T$
is represented by a stable map with image in 
$X^1$, and $e_T(N_F^\vir)$ is determined
by the $T$-equivariant formal scheme $\hat{X}^1$. 
We summarize this paragraph in Figure 5.

\medskip

\begin{figure}[h]  \label{Figure5}
\psfrag{planar trivalent graph}{\small $\begin{array}{c}
\textup{locally planar}\\
\textup{trivalent graph }\Gamma_X
\end{array}$}
\psfrag{formal neighborhood}{$\small \begin{array}{r}
T\textup{-equivariant formal}\\
\textup{neighborhood }\hat{X}^1\\
\textup{of }X^1\textup{ in }X
\end{array}$}
\psfrag{GW invariants}{$\small \begin{array}{r}
\textup{Gromov-Witten}\\
\textup{invariants }N^X_{g,\vd}
\end{array}$}
\psfrag{E1}{\small $E_1$}
\psfrag{E2}{\small $E_2$}
\psfrag{E3}{\small $E_3$}
\psfrag{E4}{\small $E_4$}
\psfrag{C0}{\small $C_0$}
\begin{center}
\includegraphics[scale=0.5]{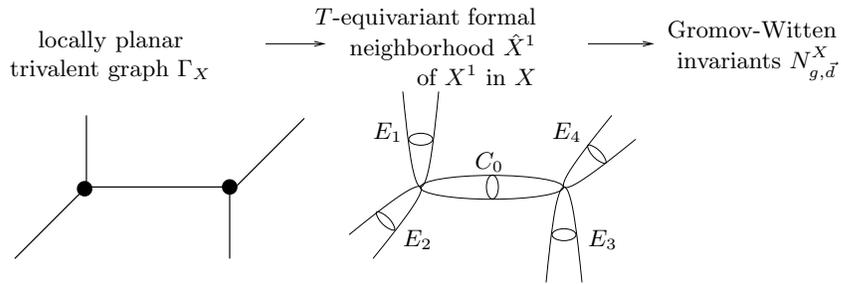}
\caption{Gromov-Witten invariants from a locally trivalent graph}
\end{center}
\end{figure}

\subsection{Formal Toric Calabi-Yau (FTCY) graphs}
\label{sec:FTCY}

In order to develop the mathematical 
theory of the topological vertex,
we need to generalize local planar trivalent
graphs associated to toric Calabi-Yau 3-folds. The generalization will be 
FTCY (formal toric Calabi-Yau) graphs.
The reverse procedure $\Gamma_X \to \hat{X}^1$ will
be generalized to construction
of relative FTCY 3-folds. We define
relative FTCY invariants by localization.
These three ingredients are summarized
in  Figure 6, which is generalization Figure 5.

\begin{figure}[h] \label{Figure6}
\begin{center}
\psfrag{FTCY}{\small FTCY graph $\Gamma$}
\psfrag{3fold}{\small $\begin{array}{c}
\textup{relative FTCY 3-fold}\\
Y_\Gamma^{\mathrm{rel}}
=(\hat{Y},\hat{D})\\
K_{\hat Y}+\hat{D}=0
\end{array}$}
\psfrag{GW}{\small $\begin{array}{c}
\textup{formal relative}\\
\textup{Gromov-Witten}\\
\textup{invariants }
F^\Gamma_{g,\vec{d},\vec{\mu}}
\end{array}$}
\psfrag{e0}{\small$w_0$}
\psfrag{-e0}{\small$-w_0$}
\psfrag{e1}{\small$w_1$}
\psfrag{e2}{\small$w_2$}
\psfrag{e3}{\small$w_3$}
\psfrag{e4}{\small$w_4$}
\psfrag{f1}{\small$f_1$}
\psfrag{-f1}{\small$-f_1$}
\psfrag{f2}{\small$f_2$}
\psfrag{-f2}{\small$-f_2$}
\psfrag{f3}{\small$f_3$}
\psfrag{-f3}{\small$-f_3$}
\psfrag{f4}{\small$f_4$}
\psfrag{-f4}{\small$-f_4$}
\psfrag{f1f2}{\small$\begin{array}{c}
f_1=w_2-n_1 w_1\\
f_2=w_0-n_2 w_2
\end{array}$}

\psfrag{f3f4}{\small$\begin{array}{c}
f_3=w_4-n_3 w_3\\
f_4=-w_0-n_4 w_4
\end{array}$}
\psfrag{D1}{\small$D_1$}
\psfrag{D2}{\small$D_2$}
\psfrag{D3}{\small$D_3$}
\psfrag{D4}{\small$D_4$}

\psfrag{C0}{\small$C_0$}
\psfrag{C1}{\small$C_1$}
\psfrag{C2}{\small$C_2$}
\psfrag{C3}{\small$C_3$}
\psfrag{C4}{\small$C_4$}
\psfrag{N=O(n)+O(-n-1)}{\small$\begin{array}{c}N_{C_i/\hat{Y}}=\cO_{\bP^1}(n_i)\oplus \cO_{\bP^1}(-n_i-1)\\
 i=1,2,3,4 \end{array}$}
\includegraphics[scale=0.5]{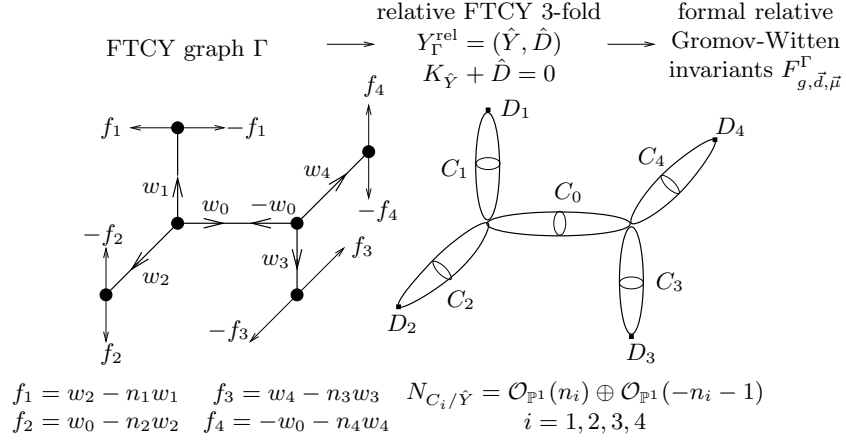}
\caption{formal relative Gromov-Witten invariants from a 
FTCY (formal toric Calabi-Yau) graph}
\end{center}
\end{figure}

As an example, we start with the graph of $\cO_{\bP^1}(-1)\oplus \cO_{\bP^1}(-1)$
on the left hand side of Figure 5. We first compactify the rays by adding a univalent
vertex at the end of each ray. We obtain a graph with 2 trivalent vertices
and 4 univalent vertices, as on the left hand side of of Figure 6. 
This corresponds to adding a point
of infinity to each $E_i\cong \bC$ in Figure 5 
so that it becomes $C_i\cong \bP^1$ in Figure 6.
We obtain a configuration of
five $\bP^1$'s, as on the right hand side of Figure 6. 
This configuration sits in a 3-dimensional formal scheme $\hat{Y}$, 
and the 4 points added are the intersection points
with the 4 connected components $D_1, D_2, D_3, D_4$ 
of the relative divisor $\hat{D}\subset \hat{Y}$. 
By the relative Calabi-Yau condition 
$K_{\hat{Y}}+\hat{D}=0$, the degrees of
$N_{C_i/\hat{Y}}$ are $-1$ for $i=1,2,3,4$ and
$-2$ for $i=0$. We introduce, at each univalent vertex,
a {\em framing} vector which determines
$n_i$ together with the $T$-action on
$N_{C_i/\hat{Y}}$ for $i=1,2,3,4$ (see Figure 6).

The formal relative Gromov-Witten invariants $F^\Gamma_{g,\vd,\vmu}$ count morphisms
$u: C\to \hat{Y}$, where $C$ is a curve of arithmetic genus $g$, 
$u_*[C]=\vd=d[C_0]+\sum_{i=1}^4 d_i[C_i]$, and the ramification patterns of $u$ along $\hat{D}$
are $\vmu=(\mu^1,\ldots,\mu^4)$. In our example, $d_i=|\mu^i|$, where
$|\mu^i|= \mu^i_1+ \mu^i_2 +\cdots$ is the {\em size} of the partition $\mu^i$.

A priori the formal relative Gromov-Witten invariants  $F^\Gamma_{g,\vd,\vmu}$ depend on 
equivariant parameters. It was proved in  \cite{LLLZ} that 
$F^\Gamma_{g,\vd,\vmu}$ are rational numbers independent of equivariant parameters, 
so they are {\em topological} invariants instead of {\em equivariant} invariants. It is crucial
to pick the subtorus $T$ because $F^\Gamma_{g,\vd,\vmu}$ would depend on 
equivariant parameters if we used the rank $3$ torus $(\bC^*)^3$.

For our purpose, we would like to introduce normal crossing singularities
of the form 
\begin{equation}\label{eqn:node}
\{(x,y)\in \bC^2\mid  xy=0 \} \times \bC^2. 
\end{equation}
For example, we degenerate
the smooth relative FTCY 3-fold $\hat{Y}$ in Figure 6,
so that $C_0$ degenerates to $C_0'$ and $C_0''$ intersecting
at a node $p$,  as on the right of  Figure 7. This degeneration corresponds to 
inserting a bivalent vertex $v$, as shown on the left of Figure 7.
The normal bundle $N_{C_0/\hat{Y}}\cong \cO_{\bP^1}(n_0)\oplus \cO_{\bP^1}(-n_0-2)$
degenerates into two degree $-1$, rank $2$ bundles over $C_0'$ and $C_0''$: 
\begin{equation}\label{eqn:split}
N_{C_0'/\hat{Y}'}\cong \cO_{\bP^1}(a)\oplus \cO_{\bP^1}(-a-1),\quad
N_{C_0''/\hat{Y}''}\cong \cO_{\bP^1}(b)\oplus \cO_{\bP^1}(-b-1),\quad a+b=n_0,
\end{equation}
where $\hat{Y}'$ and $\hat{Y}''$ are the two irreducible components
of the singular relative FTCY 3-fold. To specify the splitting
\eqref{eqn:split} $T$-equivariantly, we introduce
a framing vector $f$ at the bivalent vertex $v$.  

\begin{figure}[h] \label{Figure7}
\begin{center}
\psfrag{f}{\small $f$}
\psfrag{-f}{\small $-f$}
\psfrag{v}{\small $v$}
\psfrag{p}{\small $p$}
\psfrag{C0}{\small $C_0'$}
\psfrag{C0p}{\small $C_0''$}
\psfrag{C1}{\small $C_1$}
\psfrag{C2}{\small $C_2$}
\psfrag{C3}{\small $C_3$}
\psfrag{C4}{\small $C_4$}
\includegraphics[scale=0.5]{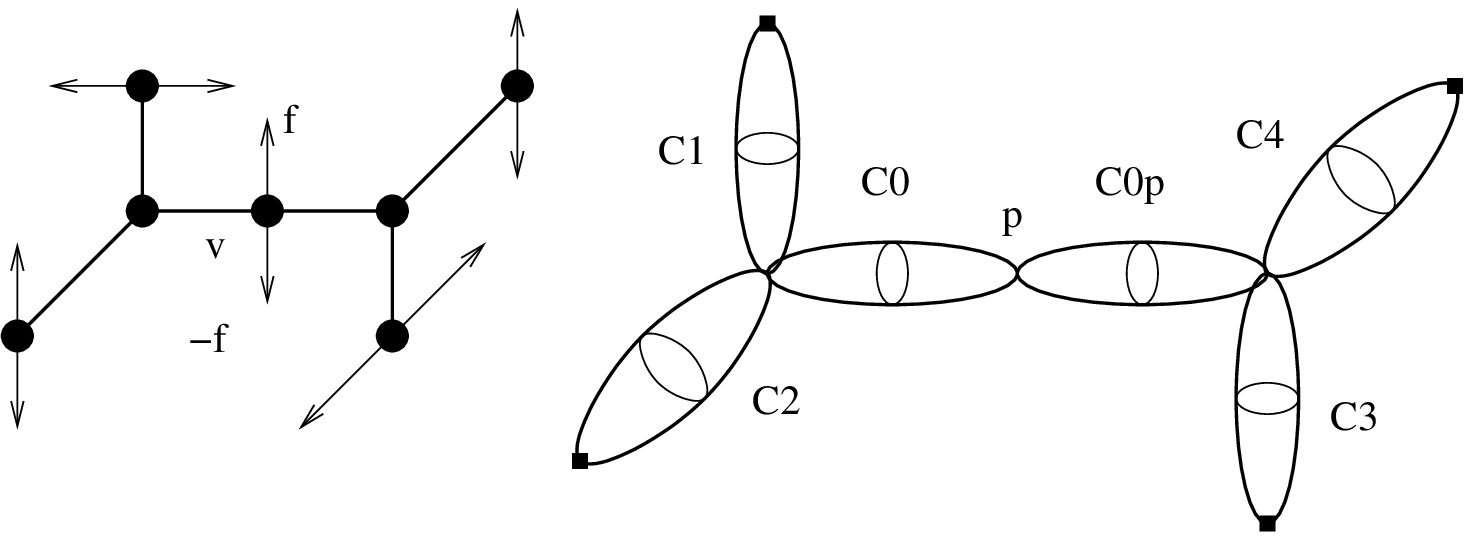}
\end{center}
\caption{ }
\end{figure}

\subsection{Degeneration formula}
\label{sec:degenerate}
Given a FTCY graph $\Gamma$ with a bivalent vertex
(see Figure 7), so that the relative FTCY 3-fold
$\hat{Y}^{\mathrm{rel}}_\Gamma$ has a normal
crossing singularity, we may either deform
or resolve this singularity.
Figure 8 shows the corresponding operation
on the graph $\Gamma$: $\Gamma_0$ is 
the smoothing of $\Gamma$, and $\Gamma_1 \cup \Gamma_2$
is the resolution of $\Gamma$.
The relative FTCY 3-folds 
$$
\hat{Y}^\mathrm{rel}_{\Gamma_0},\quad
\hat{Y}^\mathrm{rel}_{\Gamma_1\cup \Gamma_2}
=\hat{Y}^\mathrm{rel}_{\Gamma_1}\cup
\hat{Y}^\mathrm{rel}_{\Gamma_2} \textup{ (disjoint union)}
$$
are smooth. The degeneration formula relates
the relative formal Gromov-Witten invariants of 
the smoothing $\hat{Y}^\mathrm{rel}_{\Gamma_0}$ to those of  
the resolution $\hat{Y}^\mathrm{rel}_{\Gamma_1}\cup \hat{Y}^\mathrm{rel}_{\Gamma_2}$.

\begin{figure}[h]\label{Figure8}
\begin{center}
\psfrag{m1}{\small$\mu^1$}
\psfrag{m2}{\small$\mu^2$}
\psfrag{m3}{\small$\mu^3$}
\psfrag{m4}{\small$\mu^4$}
\psfrag{G0}{\small$\Gamma_0$}
\psfrag{G1}{\small$\Gamma_1$}
\psfrag{G2}{\small$\Gamma_2$}
\psfrag{G}{\small$\Gamma$ }
\psfrag{n}{\small$\nu$}
\psfrag{nd}{\small$\nu\vdash d$}
\psfrag{d}{\small$d$}
\psfrag{v}{\small$v$}
\includegraphics[scale=0.5]{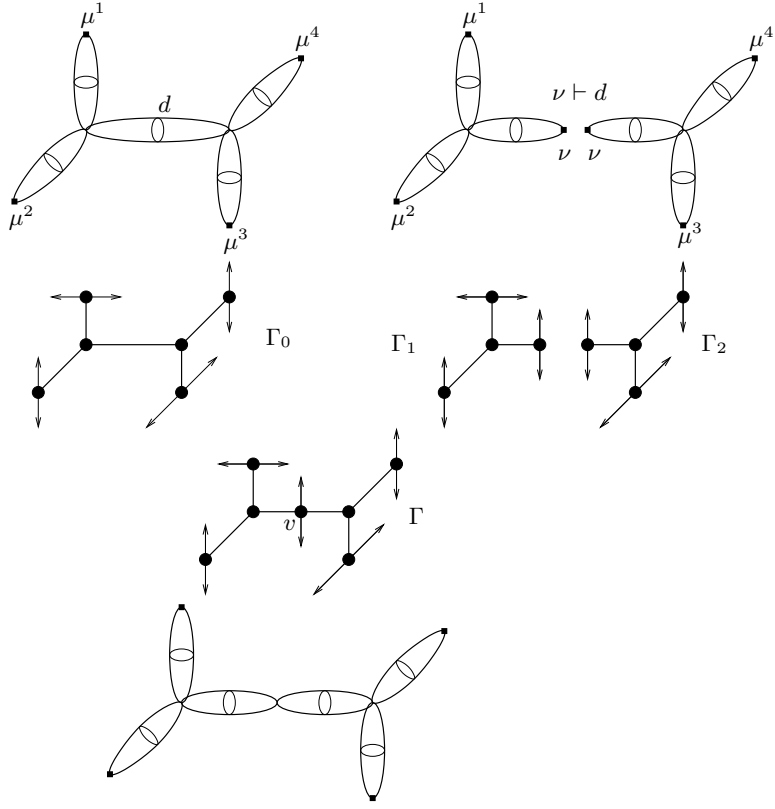}
\end{center}
\caption{Smoothing and Resolution}
\end{figure}

Geometrically, given a relative stable map 
$\hat{Y}^\mathrm{rel}_{\Gamma_0}$ with 
$$
\vd=d[C_0]+ d_1[C_1]+\cdots+ d_4[C_4],\quad
\vmu=(\mu^1,\mu^2,\mu^3,\mu^4),
$$
we degenerate and resolve to obtain
a pair of maps to $\hat{Y}^\mathrm{rel}_{\Gamma_1}$
and to $\hat{Y}^\mathrm{rel}_{\Gamma_2}$ with
degrees
$$
d[C_0']+ d_1[C_1]+ d_2[C_2],\quad d[C_0''] + d_3[C_1]+ d_4[C_4],
$$
and ramification patterns 
$$
(\mu^1,\mu^2,\nu),\quad (\nu,\mu^3,\mu^4).
$$  
(See Figure 8.) In this example, 
$d_i= |\mu^i|$ for $i=1,2,3,4$, so
we may suppress $d_1,d_2,d_3,d_4$.

Conversely, given a pair of maps to 
$\hat{Y}^\mathrm{rel}_{\Gamma_1}$ and
$\hat{Y}^\mathrm{rel}_{\Gamma_2}$, with
ramifications matching in the middle,
we may glue and smooth  to obtain
a map to $\hat{Y}^\mathrm{rel}_{\Gamma_0}$.
The degeneration formula says exactly this:
$$
F^{\bullet\Gamma_0}_{\chi,d,\vec{\mu}}
=\sum_{\tiny\begin{array}{c}
\chi_1,\chi_2,\nu\vdash d\\
\chi_1+\chi_2-2\ell(\nu)=\chi\end{array}}
F^{\bullet\Gamma_1}_{\chi_1,\mu^1,\mu^2,\nu}z_\nu
F^{\bullet\Gamma_2}_{\chi_2, \nu,\mu^3,\mu^4}
$$
where $z_\nu$ is some combinatorial factor.
Here we consider $F^\bu_{\chi, \ldots}$ which counts
maps from possibly disconnected curves $C$
with $\chi= 2\chi(\cO_C)$ 
instead of $F_{g,\ldots}$ which counts maps from
connected curves of genus $g$,  because the degeneration
formula of $F^\bu_{\chi,\ldots}$ is neater than
that of $F_{g,\ldots}$.

\subsection{Topological vertex}
\label{sec:vertex}

By degeneration and resolution,
it remains to compute the formal
invariants of the graph on the left hand side of Figure 9.
The graph depends on three integers
$n_1, n_2, n_3$ and two weights
$w_1, w_2$; note that $w_3=-w_1-w_2$
and the framings $f_i$'s are determined by 
$n_i$'s and $w_i$'s.

\begin{figure}[h]\label{Figure9}
\begin{center}
\psfrag{w1}{\small$w_1$}
\psfrag{w2}{\small$w_2$}
\psfrag{w3}{\small$w_3$}
\psfrag{f1}{\small$f_1$}
\psfrag{f2}{\small$f_2$}
\psfrag{f3}{\small$f_3$}
\psfrag{-f1}{\small$-f_1$}
\psfrag{-f2}{\small$-f_2$}
\psfrag{-f3}{\small$-f_3$}
\psfrag{m1}{\small$\mu^1$}
\psfrag{m2}{\small$\mu^2$}
\psfrag{m3}{\small$\mu^3$}
\psfrag{C1}{$C_1$}
\psfrag{C2}{$C_2$}
\psfrag{C3}{$C_3$}
\psfrag{f1=w2-nw1}{\small$f_i=w_{i+1}-n_i w_i$}
\psfrag{N=O(n)+O(-n-1)}{\small
$N_{C_i/\hat{Y}}=\cO(n_i)\oplus\cO(-n_i-1)$}
\psfrag{gamma}{$\Gamma_{w_1,w_2,\mathbf{n}}$}
\includegraphics[scale=0.6]{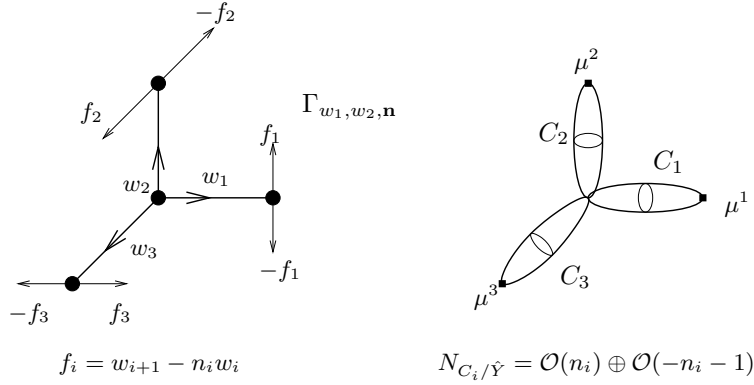}
\end{center}
\caption{The topological vertex}
\end{figure}

The invariants depend on the genus $g$ of the domain
and three partitions $\mu^1$, $\mu^2$, $\mu^3$
corresponding to the ramification patterns along
$D^1$, $D^2$, $D^3$. Note that in this case
the degrees $\vd$ is determined by ramification
patterns $\vmu$. 

Define 
$$
F^{\bullet}_{\chi,\vmu}(\vn):=
F^{\bullet\Gamma_{w_1,w_2,\vn}}_{\chi,\vmu}
$$
where $\vmu=(\mu^1,\mu^2,\mu^3)$, 
$\vn=(n_1,n_2,n_3)$.
A priori it depends on both the topological
data $n_i$ and the equivariant data $w_1,w_2$.
In \cite{LLLZ}, the authors proved that it is indeed {\em topological}:
it is a rational number depending on $n_i$ but
not on $w_i$.

$F^\bu_{\chi,\vmu}(\vn)$ can be viewed
as local relative Gromov-Witten invariants
of a configuration of three $\bP^1$'s 
embedded in a relative 
Calabi-Yau 3-fold such that
the formal neighborhood is the relative FTCY
3-fold defined by the graph of a topological vertex.
We expect the following two counting problems to be 
equivalent:
\begin{enumerate}
\item[(i)] {\em Relative topological vertex.} Counting multiple covers of
a configuration of three spheres $C_1\cup C_2\cup C_3$ (Figure 9) 
embedded in a Calabi-Yau 3-fold $Y$ relative to three divisors 
$D_1, D_2, D_3$  with ramification patterns $\mu^1,\mu^2,\mu^3$.
\item[(ii)] {\em Open topological vertex.} Counting
multiple covers of a configuration of three discs 
$D_1\cup D_2\cup D_3$ (Figure 1)
embedded in a Calabi-Yau 3-fold $Y$ relative to three
Lagrangian submanifolds $L_1,L_2, L_3$ with winding
numbers $\mu^1$, $\mu^2$, $\mu^3$.
\end{enumerate}

It is interesting to compare the above (i) and (ii) with
the classical Hurwitz problem:
\begin{enumerate}
\item[(i)'] Counting ramified covers of the
sphere by compact Riemann surfaces, with
prescribed ramification pattern $\mu$ over
$\infty$.

\item[(ii)'] Counting ramified covers of
the disk by bordered Riemann surfaces,
with prescribed winding  numbers $\mu$.
\end{enumerate}
The above two counting problems (i)' and (ii)'
are equivalent; they both give rise to Hurwitz numbers.

We now introduce some notation.
Given a partition $\mu=(\mu_1 \geq \cdots \geq \mu_n >0)$,
let $\ell(\mu)=n$ be the {\em length} of the partition.
Given a triple of partitions $\vmu=(\mu^1,\mu^2,\mu^3)$, let
$\ell(\vmu)=\ell(\mu^1)+\ell(\mu^2)+\ell(\mu^3)$, and 
define a generating function
$$
F^\bu_{\vmu}(\lam;\vn)=
\sum_{\chi}\lam^{-\chi+\ell(\vmu)}F^\bu_{\chi,\vmu}(\vn).
$$

\subsection{Localization} 
For later convenience, we consider a slightly modified generating function
$$
\tF_{\vmu}^\bu(\lam;\vn)
=(-1)^{\sum_{i=1}^3(n_i-1)|\mu^i|}
\sqrt{-1}^{\ell(\vmu)}F^\bu_{\vmu}(\lam;\vn).
$$
This is a generating function of relative Gromov-Witten invariants in the {\em winding basis}.
By localization calculations, these invariants can be expressed in terms
of {\em three-partition Hodge integrals} and {\em double Hurwitz numbers}, which
we define now.

The {\em three-partition Hodge integrals} are defined to be
\begin{eqnarray*}
G_{g,\vmu}(w_1,w_2,w_3)
&=&\frac{(-\sqrt{-1})^{\ell(\vmu)}}{|\Aut(\vmu)|}
\prod_{i=1}^3\prod_{j=1}^{\ell(\mu^i)}
\frac{\prod_{a=1}^{\nu_j-1}(w_{i+1}\mu^i_j+aw_i)}
{(\mu^i_j-1)!w_i^{\mu^i_j-1}}\\
&&\cdot\int_{\Mbar_{g,\ell(\mu)}}\prod_{i=1}^3
\frac{\Lambda_g^\vee(w_i)w_i^{\ell(\vmu)-1}}
{\prod_{j=1}^{\ell(\mu^i)}(w_i(w_i-\mu^i_j\psi_{l_i+j}))}
\end{eqnarray*}
where 
$$
\ell_1=0,\ \ \ell_2=\ell(\mu^1),\ \ \ell_3=\ell(\mu^1)+\ell(\mu^2)
$$
$$
\Lambda^\vee_g(u)=u^g-\lambda_1 u^{g-1}+\cdots+(-1)^g\lam_g
$$
Define a generating function 
$$
G_{\vmu}(\lambda;w_1,w_2, w_3)=\sum \lambda^{2g-2+\ell(\vmu)}
G_{g,\vmu}(\lambda; w_1, w_2, w_3).
$$
Let $G^\bu_{\vmu}(\lambda; w_1,w_2,w_3)$ be the disconnected version
of $G_{\vmu}(\lambda;w_1,w_2,w_3)$.

Let $H^\bu_{\chi,\nu,\mu}$ be disconnected double Hurwitz numbers
which count possibly disconnected 
covers of $\bP^1$ with ramification patterns $\nu$ and $\mu$
over $0$ and $\infty$. We define a generating function
$$
\Phi^\bu_{\nu,\mu}(\lambda)=
\sum_\chi H^\bu_{\chi,\nu,\mu}
\frac{\lambda^{-\chi+\ell(\nu)+\ell(\mu)}}{(-\chi+\ell(\nu)+\ell(\mu))!}.
$$

Localization calculations yield the following expression:
\begin{equation}\label{eqn:FGPhi}
\tF_{\vmu}^\bu(\lam;\vn)=
\sum_{ |\nu^i|=|\mu^i| }G^\bu_{\vnu}(\lam;w_1,w_2,w_3)
\cdot\prod_{i=1}^3 z_{\nu^i}\Phi^\bu_{\nu^i,\mu^i}
\Bigl(\sqrt{-1}\lambda\bigl(n_i-\frac{w_{i+1}}{w_i}\bigr)\Bigr ).
\end{equation}

\subsection{Framing dependence}

Let $d=|\nu|=|\mu|$. 
Recall that for each partition $\mu$ of $d$, there is an associated
irreducible character $\chi_\mu$ of the symmetric group $S_d$
and an associated conjugacy class $C_\mu$ of $S_d$.
By Burnside formula,
\begin{equation}\label{eqn:burnside}
\Phi^\bu_{\nu,\mu}(\lambda)=
\sum_{\sigma \vdash d}e^{\kappa_\sigma\lambda/2}
\frac{\chi_\sigma(C_\nu)}{z_\nu} \frac{\chi_\sigma(C_\mu)}{z_\mu}.
\end{equation}
Equations \eqref{eqn:FGPhi} and \eqref{eqn:burnside} imply 
the following.

\begin{pro}[framing dependence in the winding basis]
$$
\tF_{\vmu}^\bu(\lambda;\vn)
=\sum_{|\nu^i|=|\mu^i|}\tF_{\vnu}^\bu(\lambda;\mathbf{0})
\prod_{i=1}^3 z_{\nu^i}\Phi^\bu_{\nu^i,\mu^i}(\sqrt{-1} \lambda n_i).
$$
\end{pro}

We introduce a generating function
$$
\tC_\vmu(\lam;\vn)=\sum_{|\nu^i|=|\mu^i|}\tF^\bu_\vnu(\lam;\vn)
\prod_{i=1}^3\chi_{\mu^i}(C_{\nu^i}).
$$
These are relative Gromov-Witten invariants in the
{\em representation basis}.
The framing dependence in the representation
basis is simple:
\begin{pro}[framing dependence in the representation basis]
$$
\tC_\vmu(\lam;\vn)=q^{\sum \kappa_{\mu^i}n_i/2}\tC_{\vmu}(\lam;\mathbf{0}),
\quad q=e^{\sqrt{-1}\lam}.
$$
\end{pro}

Note that O3 in Section \ref{sec:AKMV} implies 
$$
C_{\vmu}(\lam;\vn)=q^{\sum\kappa_{\mu^i}n_i/2} C_{\vmu}(\lam;\mathbf{0}).
$$

Up to now, we have defined invariants $\tilde{C}_\vmu$ which
have the same gluing formula and framing dependence  as $C_\vmu$ do.

\subsection{Combinatorial expression} \label{sec:tW}
We have
\begin{lm}
\begin{equation}\label{eqn:Gg-reduce}
\begin{aligned}
G_{g,\vmu}(\lambda;1,1,-2)
=& (-1)^{|\mu^1|-\ell(\mu^1)}
\frac{z_{\mu^1\cup \mu^2}}{z_{\mu^1}\cdot z_{\mu^2}}
G_{g,\emptyset,\mu^1 \cup \mu^2, \mu^3}(\lambda;1,1,-2)\\
& +\delta_{g0}\sum_{m\geq 1} \delta_{\mu^1 (m)} \delta_{\mu^2\emptyset}
\delta_{\mu^3(2m)}\frac{(-1)^{m-1}}{m}
\end{aligned}
\end{equation}
\end{lm}

A formula of the two-partition Hodge 
integrals $G_{g, \emptyset, \mu, \nu}$
in terms of $\cW_{\mu \nu}$ (the
Chern-Simons invariants of the Hopf link)
was derived in \cite{LLZ2}. 
$\cW_{\mu \nu}$ can be expressed in
terms of the skew functions
$s_{\mu/\lambda}$ (see \cite{Ok-Re-Va}):
$$
\cW_{\mu\nu}(q)=q^{(\kappa_\mu+\kappa_\nu)/2}\sum_\lambda
s_{\mu^t/\lambda}(q^{-\frac{1}{2} }, q^{-\frac{3}{2}},\ldots)
s_{\nu^t/\lambda}(q^{-\frac{1}{2}} , q^{-\frac{3}{2}},\ldots).
$$
This allows one to evaluate the relative Gromov-Witten
invariants $\tC_{\vmu}$ of the topological vertex  in terms 
of $\cW_{\mu \nu}(q)$. More precisely, let $\cW_\mu(q)=\cW_{\mu,\emptyset}(q)$,
and let $c^\mu_{\eta\rho}$ be the Littlewood-Richardson
coefficients. Then
\begin{enumerate}
\item[R3.] $\tC_{\vmu}(\lam;\mathbf{0})=\tilde{\cW}_{\vmu}(q)$, where
\begin{eqnarray*}
\tilde{\cW}_{\rho_1,\rho_2,\rho_3}(q)
&=&q^{-(\kappa_{\rho_1}-2\kappa_{\rho_2}-\frac{1}{2}\kappa_{\rho_3})/2}
 \sum c^{\nu^+}_{(\nu^1)^t\rho^2}
c^{\rho^1}_{(\eta^1)^t\nu^1}
c^{\rho^3}_{\eta^3\nu^3}\\
&&\cdot q^{(-2\kappa_{\nu^+}-\frac{\kappa_{\nu^3} }{2})/2}\cW_{\nu^+,\nu^3}(q)
\frac{1}{z_\mu}\chi_{\eta^1}(\mu)\chi_{\eta^3}(2\mu)
\end{eqnarray*}
\end{enumerate}

\subsection{Applications}
We now have an explicit formula of all
formal relative Gromov-Witten invariants of relative FTCY
3-folds (and in particular, Gromov-Witten invariants
of all toric Calabi-Yau 3-folds) in terms of $\tilde{\cW}_\vmu$. 
This formula has computational and theoretical applications:
\begin{enumerate}
\item[A1.] It is significantly more efficient than the 
traditional algorithm described in Section \ref{sec:tradition}.
Given a toric Calabi-Yau 3-fold, 
the generating function of its Gromov-Witten invariants 
in all genera in a fixed degree is an infinite series.
By this formula, this infinite series is equal
to a finite sum in terms of symmetric functions.

\item[A2.] This formula can be used to prove structural
theorems of Gromov-Witten invariants. For example, P. Peng used 
this formula to prove the  Gopakumar-Vafa conjecture \cite{Go-Va2} for
toric Calabi-Yau 3-folds \cite{Pe}.
\item[A3.] One can use this formula to verify
enumerative predictions by other string dualities, for
example the geometric engineering (see \cite{Li-Liu-Z}).
\end{enumerate}

\subsection{Comparison}\label{sec:W}
The physical theory of the topological vertex
predicts the following formula for $C_{\vmu}(q)$:
\begin{enumerate}
\item[O3.] $C_{\vmu}(\lam;\mathbf{0})=\cW_{\vmu}(q)$, where
$$
\cW_{\mu^1,\mu^2,\mu^3}(q)
=q^{(\kappa_{\mu^2}+\kappa_{\mu^3})/2}
\sum c^{\mu^1}_{\eta\rho^1} c^{(\mu^3)^t}_{\eta (\rho^3)^t}
\frac{\cW_{(\mu^2)^t \rho^1}(q)\cW_{\mu^2(\rho^3)^t}(q) }{\cW_{\mu^2}(q)}.
$$
\end{enumerate}

The equivalence of the physical theory and the mathematical
theory of the topological vertex boils down to the
following identity of classical symmetric functions:
\begin{equation}\label{eqn:WtW}
\tilde{\cW}_\vmu(q)=\cW_\vmu(q).
\end{equation}
Equation \eqref{eqn:WtW} follows from the results in \cite{MOOP} (see Section \ref{sec:GWDT}
below).

\section{GW/DT Correspondence and the Topological Vertex}
\label{sec:GWDT}
Maulik, Nekrasov, Okounkov, and Pandharipande conjectured a
correspondence between the GW (Gromov-Witten)
and DT (Donaldson-Thomas) theories for any non-singular
projective 3-fold \cite{MNOP1, MNOP2}.
This correspondence can also be formulated for
certain noncompact 3-folds in the presence
of a torus action; the correspondence for toric Calabi-Yau 3-folds is
equivalent to the algorithm of the topological vertex \cite{MNOP1, Ok-Re-Va}.
For non-Calabi-Yau toric 3-folds the building block is the equivariant vertex
(see \cite{MNOP1, MNOP2, Pa-Th1, Pa-Th2}) which depends on equivariant 
parameters. The GW/DT correspondence for all toric 3-folds
has been proved in a recent work by Maulik-Oblomkov-Okounkov-Pandharipande \cite{MOOP}.

\bigskip

\subsection*{Acknowledgments} I am very grateful to Professor Shing-Tung Yau, who
has educated me and has been offering me tremendous help and support over the years.
As my thesis advisor, Professor Yau encouraged me to study the mathematical aspect of the 
large $N$ duality and provided me opportunities to learn its newest developments from leading experts.
It is my great pleasure to dedicate this article to Professor Yau on the occasion of his 59th birthday. 

I would like to thank my collaborators Jun Li, Kefeng Liu, and Jian Zhou for the fruitful 
and pleasant collaboration.

\end{document}